\documentclass[11pt,a4paper]{article}

\title{On quantum topology, hypergraphs and flag vectors}

\author{Jonathan Fine\relax
\thanks{203 Coldhams Lane, Cambridge, CB1 3HY, England.
\hfil\break
\qquad
E-mail: \texttt{j.fine@pmms.cam.ac.uk}}%
}
\date{1 August 1997}

\newcommand\calG{{\cal G}}
\newcommand\calM{{\cal M}}
\newcommand\calK{{\cal K}}
\newcommand\calX{{\cal X}}
\newcommand\ftilde{\widetilde{f}}
\newcommand\Ghat{\widehat{G}}
\newcommand\bdot{b^{\textstyle\cdot}}
\newcommand\bfR{{\bf R}}
\newcommand\bibdash{\vrule width 2.5pc depth 0pt height 0.4pt \relax}

\begin{document}
\maketitle

\begin{abstract}\noindent
Each rule $f$ that assigns a vector $f(G)$ to an $(n+1)$-graph $G$
determines a class (or property) of $n$-manifold invariants.  An invariant
$v=v(M)$ is in this class if, for any triangulated manifold $|G|=M$, one
has that $v(M)$ is a linear function of $f(G)$.  This paper defines a flag
vector $f(G)$ for $i$-graphs, whose associated invariants might be
quantum, and which is of interest in its own right. The definition (via
the concept of shelling, and a `disjoint pair of optional cells' rule for
the link) seems to apply to any finite combinatorial object, and so to any
compact topological object that can be triangulated.  It also applies to
finite groups.
\end{abstract}

\section {Introduction}

The purpose of this paper is to formulate some combinatorial questions
related to the theory of quantum manifold invariants, and make progress
towards their solution.  Briefly, the problem is this.  Suppose $v =
v(M)$ assigns a number, or more generally a vector, to every compact
topological manifold $M$ of some fixed dimension $n$.  In particular, if
$T$ is a suitable simplicial complex then its geometric realisation $|T|$
will be such a manifold.  Thus, the quantity $v(T)=v(|T|)$ is defined for
certain simplicial complexes $T$.  The problem is to extend $v$ in a
useful and instructive way, so that it applies to more general simplicial
complexes.

Quantum topology is not yet fully developed.  For this reason, it may be
best in the above to treat the topological invariant $v$ as an unknown. 
Instead, suppose that some rule $f=f(T)$ that assigns a vector to each
simplicial complex $T$  has been provided.  This rule determines a class
(or property) of manifold invariants, in the following way.

Each manifold $M$ can be triangulated in many ways.  Each triangulation is
a solution to the equation $|T|=M$.  For certain linear functions
$v=v(T)=v(f(T))$ of the $f$-vector of $T$, the value of $v(T)$ will depend
only on $M=|T|$, and not on the triangulation chosen.  In other words,
certain linear functions $v$ of $f$ are, on triangulated manifolds,
topological invariants.

Each rule $f=f(T)$ that assigns a $f$-vector to a simplicial complex $T$
thus determines a class of topological invariants of manifolds, namely
those that are linear functions of $f$.  The problem is to formulate a
rule that will produce the presently loosely defined class known as the
quantum invariants.  Note that it can be far easier to define a whole
class of invariants in this way, than it is to produce a single invariant
that satisfies the membership property.  This paper will suggest a
solution, but first some general remarks are in order.

Topologists, under the influence of homotopy and homology theory, have
developed the concept of a simplicial complex.  For the present purpose,
however, the concept of a hypergraph is more appropriate.  Suppose that
some finite set $V$ of vertices is given.  A \emph{simple $i$-uniform
hypergraph on $V$}, or $i$-graph for short \cite{bib:Duchet}, is simply a
collection $C$ of $i$-element subsets of the vertex set $V$.  Each element
of $C$ will be called a \emph{cell} (or edge) of the hypergraph.  (Thus, a
$2$-graph is an ordinary graph, while a triangulation of the $n$-manifold
$M$ produces an $(n+1)$-graph.  The cells are the vertex sets of the
triangulation.  A $1$-graph is a collection of $1$-element subsets, or
what amounts to much the same thing, a subset of $V$.  A cell of a
$1$-graph is then a vertex, that has been selected to appear in $C$.) In
\S2 the flag vector $fG$ of an $i$-graph will be defined.  It may be that
it corresponds to quantum topology.  Some care is taken to justify the
definition given.  To do this, some arguments to appear in \cite{bib:GFP}
will be summarized.

For this paper it is not so much the flag vector as the concept of a
linear function of the flag vector that is of interest.  This can be
expressed in the following way.  Let $\calG$ be the module of all formal
sums of $i$-graphs on some vertex set $V$.  (Sometimes, a formal sum of
$i$-graphs will be called an $i$-graph.  The prefix `actual' will be used
to recover the original concept. A formal sum of two graphs is not the
same as their disjoint union.)   Now define the \emph{nullspace} $\calG_0$
to consist of all formal sums
\[
    G = \lambda_1 G_1 + \dots + \lambda_r G_r
\]
in $\calG$, such that
\[
    f(G) = \lambda_1 f(G_1) + \dots + \lambda_r f(G_r)
\]
is zero.  A rule $v=v(G)$ that assigns a vector to each $i$-graph is then
a linear function of the flag vector $fG$ just in case
\[
    v(G) = \lambda_1 v(G_1) + \dots + \lambda_r v(G_r)
\]
is zero for every $G$ in $\calG_0$.  The flag vector is a means of
defining $\calG_0$.  Conversely, the natural map that takes a $i$-graph
$G$ to its residue $\calG/\calG_0$ assigns a vector $fG$ to each
$i$-graph, and has $\calG_0$ as its null space.  This provides the
\emph{abstract form} of the flag vector.

To be useful, the nullspace should be neither too large or too small.  If
$\calG_0=\{0\}$ then any rule $v=v(G)$ will be a linear function of the
flag vector, and the concept ceases to have any significance.  Similarly,
if $\calG_0=\calG$ then zero is the only linear function of the flag
vector, and again the concept is without value.  The flag vector defined
in this paper produces a nullspace with several attractive properties. 
It may also be that some simple geometric operations will produce a
spanning set of nullvectors.  This is evidence, independent of the
details of quantum topology, that the choice of $\calG_0$ is correct.

Now suppose that $v=v(M)$ is a vector valued invariant of topological
manifolds (and so also a function on hypergraphs whose realisation is a
manifold).  Even without extending $v$ to all hypergraphs, one can talk of
$v$ being a linear function of the flag vector.  Inside $\calG$ form
$\calM$, the span of $n$-graphs whose realisation is a manifold.  Define
the \emph{manifold nullspace} $\calM_0$ to be $\calM\cap\calG_0$.  If $v$
vanishes on $\calM_0$, then say that it is a linear function of the flag
vector.  Such a function can of course be extended in many ways to the
space $\calG$ of all hypergraphs.

The polynomial quantum knot invariants can be computed by means of
crossing change rules.  Such rules are easy to state and apply.  To find
the rules is harder, and harder yet is to show that the same answer will
result, however a calculation is done.  When expressed in Vassiliev form
(substitute $q=e^x$ and let $v_i$ be the coefficient of $x^i$), it follows
immediately from the form of the crossing change rule that each $v_i$
vanishes on a certain subspace $\calK_{i+1}$ of the space $\calK$ of all
formal sums of knots, namely that spanned by knots with $i+1$ double
points.  (This is because
\[
    (\sqrt{q} - 1/\sqrt{q}) = x + x^3/24 + \mbox{higher order}
\]
and so each double point contributes a factor of $x$ to a product.  This
argument appears in \cite[Thms.~2 and~3]{bib:Bar-Natan}.)

One might wish for a similar theory, applying to manifolds.  For knots,
both the crossing change formulae and the definition of the nullspaces
$\calK_i$ are easy to state, if not discover.  For manifolds the situation
is not so clear.  In this context, see \cite[Prop.~4.6]{bib:Le}.  The flag
vector $fG$ of a hypergraph can be `filtered by complexity' to produce
`components' $f_iG$, and this can be used to define subspaces $\calG_i$
and $\calM_i$ of $\calG$ and $\calM$ respectively. From here one may be
able to work backwards, and find first a geometric characterisation of
$\calM_i$, and then a `crossing change' formula for some given quantum
topological invariant.  This interplay between crossing change formulae
and nullspaces seems to be particularly important.

\section {Hypergraphs}

It is now time to turn to the definition of the flag vector $fG$ of a
hypergraph $G$.  This will be done in two stages.  First, the shelling
vector $\ftilde G$ will be defined.  Next, the flag vector $fG$ is
defined, as a variant on the construction of $\ftilde G$.  Suppose that
$G$ is an $i$-graph on $N$ vertices.  Now remove the vertices from $G$,
one at a time, until none are left.  This, together with a record of the
changes that occur at each step, is a \emph{shelling} of $G$.  It is
analogous to Morse theory in differential topology.

When a vertex $v$ is removed from an $i$-graph, the cells that contain
that vertex must also be removed.  Thus, a shelling is an ordering of the
vertices of $G$, and for each vertex $v$ a record of the cells which have
$v$ as their first vertex, for the given ordering.  One can also reverse
the process, and think of a shelling as a way of building up a graph out
of nothing.

When a vertex $v$ is removed from an $i$-graph $G$, the cells that are
removed can be described via an $(i-1)$-graph on one fewer vertices. (Each
removed cell contains $v$.  It will also contain $(i-1)$ of the remaining
vertices.)  Call this $(i-1)$-graph the \emph{link} $L_v$ at the vertex
$v$ of $G$.  Locally, about $v$, the graph $G$ looks like the cone $CL_v$
on the link $L_v$.  (The \emph{cone} on a $j$-graph is formed by adding a
new vertex to the vertex set, and also adding this vertex to all the cells
of the $j$-graph.  The result is a $(j+1)$-graph.)

The \emph{shelling vector} $\ftilde G$ is a sum over the $N!$ possible
shellings of $G$.  For each shelling sequence $\sigma$ one obtains a sequence
$L_1$, $L_2$, $\ldots$, $L_N$ of links.  The link on the $l$-th vertex is
an $(i-1)$-graph on $N-l$ vertices.  Loosely speaking, the contribution
made by $\sigma$ to $\ftilde G$ is this sequence of links, as a formal and
noncommutative product.  However, each link $L_l$ is an $(i-1)$-graph
which, by induction, can be supposed to have a shelling vector $\ftilde
L_l$.

The basis for the induction are the $0$-graphs.  These are collections of
zero element subsets of the vertex set.  There is only one zero element
subset, namely the empty set, and so there are only two possible
$0$-graphs, namely $\emptyset$ and $\{\emptyset\}$.  (The first is the
empty collection of subsets, the second a non-empty collection, whose only
element is the empty set.)  Call these $0$-graphs $[a]$ and $[b]$
respectively.

Now let $[\>,\ldots,\>]$ denote a noncommutative but linear formal
product.  The equation
\[
    \ftilde G = \sum \nolimits _\sigma \>\>
        [ \ftilde L_1 , \ftilde L_2 , \ldots , \ftilde L_N ]
\]
provides a recursive definition of the shelling vector $\ftilde G$, in
terms of $[\>,\ldots,\>]$ and the basic $0$-graphs $[a]$ and $[b]$.
(During calculations, the brackets help keep track of where one is in the
recursion.  Each term in $\ftilde G$ is an expression in $a$, $b$ and
$[\>,\ldots,\>]$.)

There is an analogy between shelling vectors and Feymann path integrals. 
In both cases one considers the totality of all ways of achieving
something, in the one case a particular hypergraph, in the other some
outward form of an interaction between particles.  The final value is a
sum of the contributions made by the various ways.  In the one case the
contribution is considered as a formal entity, in the other it is a number
depending on the paths and inner interactions of the particles.

The shelling vector $\ftilde G$ is extremely large.  It may even be that
if
\[
    G = \lambda_1 G_1 + \dots + \lambda_r G_r
\]
is a formal sum of $i$-graphs, then the associated shelling vector
\[
    \ftilde (G) = \lambda_1 \ftilde G_1 + \dots + \lambda_r \ftilde G_r
\]
is zero only when $G$ itself is zero.  If true, say that the shelling
vector \emph{distinguishes formal sums} (of $i$-graphs).  It would be
useful to know whether or not this statement is true.  It is considerably
stronger than saying that if two $i$-graphs $G_1$ and $G_2$ have the same
shelling vector, then they are isomorphic.  This is the property of
\emph{distinguishing individual $i$-graphs}.

The flag vector $fG$ is obtained by using only part of the information
available in each link $L_l$.  Indeed, if $\ftilde G$ does distinguish
formal sums of $n$-graphs, some data in the link will have to be discarded
to obtain the flag vector, if not all manifold invariants are to be linear
functions of the flag vector.

Here is an example.  The shelling vector of the $1$-graph that consists of
$m$ vertices chosen to be cells out of $n$ vertices is the sum of all ways
of multiplying $m$ copies of $b$ and $n-m$ copies of $a$.  If $n=2$ then
$2aa$, $ab+ba$, and $2bb$ are the possible shelling vectors.  (Here, the
square brackets are redundant.)  Notice that the $1$-graphs on $n$
vertices have linearly independent shelling vectors.

\emph{How many} is an important question to ask of a $1$-graph (or subset
of the vertices).  For $G$ a $1$-graph, define $f'G$ to be $a+mb$, where
$m$ is the number of cells in the graph $G$.  Here $a$ and $b$ are not
quite the same symbols are were used for the shelling vector.  (The $a$
denotes the $1$-graph itself.  In a formal sum of graphs, one would like
to know how many graphs, which is the sum of the coefficients
$\lambda_i$.)

The quantity $f'G$ is a linear function of $\ftilde G$.  Its kernel is
rather interesting.  Consider $1$-graphs on 2 vertices.  One has
$0+2=1+1$, or $(0+2)-(1+1)=0$.  For shelling vectors the corresponding
expression
\begin{equation}
    2aa + 2bb - 2(ab+ba) = 2(a-b)(a-b)
\end{equation}
lies in the kernel.  More generally, the symmetric products with two
occurrences of $(a-b)$ span the kernel.  This may be important later.

The \emph{flag vector of a $1$-graph} can now be defined.  It is the sum
\begin{equation}
    fG = \sum \nolimits _\sigma \>\>
            f'L_1 \cdot f'L_2 \cdot \ldots \cdot f'L_N
\end{equation}
over all shellings of the formal product of the corresponding \emph{link
contributions} $f'L_l$, where $f'$ is as above.  In other words, as each
shelling removes a vertex, contribute an `$a$' for the vertex being
removed, and a `$b$' for each edge that is removed.  Clearly $fG$ is a
formal sum of words in $a$ and $b$, of length $N$.  It thus has $2^N$
components.

The components of $fG$ are not independent.  For example, if $A$, $B$, $C$
and $D$ are the graphs on three vertices with $0$, $1$, $2$ and $3$ edges
respectively, then the equation $fA-3fB+3fC-fD=0$ holds.  In fact
\cite{bib:GFP}, the flag vectors of $2$-graphs on $n$ vertices have
$p(n)$, the number of partitions of $n$, as the number of independent
components.  (Thus, any given linear function $v$ of graph flag vectors
$fG$ can be written in many ways, as a linear function of words in $a$ and
$b$.)

The flag vector is defined not only for graphs, but also for formal sums
of graphs.  Certain formal sums are useful in the study of flag vectors. 
Recall that an $i$-graph is a pair $(C,V)$, where $V$ is the vertex set,
and $C$ is a collection of $i$-element subsets of $V$.  An \emph{$i$-graph
$G$ with optional cells} is defined in the following way.  Let $A$ and $B$
be disjoint collections of $i$-element subsets of $V$, or in other words,
collections of potential cells.  Form the alternating sum
\[
    G = \sum \nolimits _{ C \subseteq B} \>\>
          (-1) ^{\#B - \#C} ( A \cup C , V )
\]
of the $2^{\# B}$ graphs $( A \cup C , V )$ whose cell set lies between
$A$ and $A\cup B$.  (Here, $\#B$ means the number of elements in the set
$B$.)  To denote an optional edge on a picture of a $2$-graph, use a
dotted line.  The options are to either join ($+$) or erase ($-$) the
dots.

Now suppose that a $2$-graph contains an optional cycle, or in other words
a cycle consisting of optional edges.  It then follows that its flag
vector is zero.  (For each shelling, $f'L_l$ will be zero for the first
vertex that lies on the cycle of options.) This result, applied to a
triangle of optional edges, gives the flag vector equation just recently
stated.

For $i$-graphs, the rule for the link contribution $f'L_i$ will be chosen
by elevating to a general principle a particular aspect of the case already
considered.  Let $G_0$, $G_1$ and $G_2$ denote the $1$-graphs on $2$
vertices, with $0$, $1$, and $2$ cells (selected vertices) respectively. 
The equation
\begin{equation}
\label{eqn:0+2=1+1}
    f'G_0 + f'G_2 = f'G_1 + f'G_1
\end{equation}
can be interpreted in the following way.  Label the vertices $L$ and $R$
(for left and right).  Let $G_L$ and $G_R$ be the $1$-graphs in which
$\{L\}$ and $\{R\}$ are the only cells.  Now consider the $1$-graph in
which both $L$ and $R$ are optional cells.  This is a formal sum 
\[
    G_2 - G_L - G_R + G_0
\]
whose link contribution $f'$ is, by (\ref{eqn:0+2=1+1}), zero.

Now suppose $G = G_+ - G_-$ is a $1$-graph with one optional cell. 
Clearly, when formally considered, the number of cells in $G$ (the number
in $G_+$ less that in $G_-$) is equal to one.  Similarly if
\begin{equation}
\label{eqn:G++}
    G = G_{++} - G_{+-} - G_{-+}  + G_{--}
\end{equation}
is a $1$-graph with two optional cells, then formally considered $G$ has
zero as its number of cells.  Clearly, any $m$ celled $1$-graph on $n$ vertices can be
written as $G_0 + mG_1$, plus various terms of type (\ref{eqn:G++}).  No
other value for $m$ is possible.  (Here, $G_0$ and $G_1$ are $1$-graphs on
$n$ vertices, with $0$ and $1$ cells respectively.)

As already noted, in this paper the flag vector is of interest only as a
means of defining the concept of a linear function of the flag vector. 
(When it comes to calculations and proofs, a well chosen system of
coordinates may of course be useful.)  The same applies to the link
contributions.  With this in mind, the definition of $fG$ for $2$-graphs
will be reformulated.

It is as before a sum
\begin{equation}
\label{eqn:tensor-f}
    fG = \sum \nolimits _\sigma \>\>
            f'L_1 \otimes \ldots \otimes f'L_N
\end{equation}
over all shellings.  Each link $L_l$ is a $1$-graph on $N-l$ vertices. 
Let $V_l$ be the vector space of formal sums of such graphs (up to
isomorphism), modulo the graphs with two optional cells.  Let the link
contribution $f'L_l$ be the residue of $L_l$ in $V_l$.  Each $V_l$ has
dimension $2$ (except $V_N$, for there the vertex set is empty), and the
previous value of $a+mb$ provides a system of coordinates on $V_l$.  This
new definition differs from the old by an invertible linear
transformation.

The concept of a graph with two optional cells applies not only to
$1$-graphs, but also to $i$-graphs.  In the general case, however, a
restriction will be placed on the pairs of cells that can be used.  Say
that two cells are \emph{disjoint} if they do not share a common vertex. 
It means that for any shelling, one or the other of the pair of optional
cells will appear first.  They will not both be removed at the same time. 
The options are independent in that they do not directly interact.  This
concept is trivial for $1$-graphs.

The \emph{flag vector $fG$ of an $i$-graph}, in its tensor form, is
defined by the recursive formula (\ref{eqn:tensor-f}), where now each
$f'L_l$ is the residue of $fL_l$, modulo the graphs with a disjoint pair
of optional cells.  Thus, $f'L_l$ is defined by the $(i-1)$-graph flag
vector nullspace, augmented by the disjoint pairs of optional cells.

This defines a flag vector for $i$-graphs, and thus, as described already,
a class (or property) of invariants of compact topological manifolds.

\section {Summary and Further Problems}

Much remains to be done, before a clear connection can be established
between the hypergraph flag vector on the one side and quantum topology on
the other, assuming indeed that there is such a connection.  This section
provides a summary of what is already known, and a description of some 
untouched problems.  It is an essay that outlines possible future
developments.

To begin with what is known: for ordinary or $2$-graphs, the flag vector
has attractive properties \cite{bib:GFP}.  The space spanned by the flag
vectors on $n$ vertices has dimension $p(n)$, the number of partitions of
$n$.  The nullspace of the flag vector can be given a geometric
description.  For $n=4$, the flag vector distinguishes graphs.  The flag
vectors of the $11$ distinct such graphs are the vertices of a convex
polytope in an affine $4$ ($=p(4)-1$) dimensional space. (Thus, the affine
linear functions that define the polytope are non-negative on the flag
vectors of actual graphs.)  Methods similar to those used for convex
polytopes may make similar subtle linear inequalities accessible, for
larger $n$.  Finally, if the realisation $|G|$ is a $1$-manifold ($|G|$ is
a disjoint union of polygons), the number of components in $|G|$ is a
linear function of the flag vector $fG$.  (This linear function is quite
complicated, when expressed in $a$ and $b$ form.  This is because $fG$ has
only a restricted view of the connectivity of $G$.)

For topology, quantum or otherwise, the three basic questions are these. 
First, \emph{does the shelling vector encode significant topological
information?}  If it distinguishes formal sums of $i$-graphs then the
answer is, of course, yes.  The second question is this:  \emph{does the
flag vector make this information available in a useful form?}  For
example, can the dimension of the span of the flag vectors (and, in
coordinates, the span itself) be easily described.  The third question
is:  \emph{what part of the flag vector is a topological invariant of the
realisation?}

Suppose now that there is significant information in the shelling vector. 
A major problem in using it is this.  Recall that the shelling vector
$\ftilde G$ of an $n$-graph is a formal sum of words, where each word is a
product of the shelling vectors $\ftilde L_l$ of the links $L_l$.  However
the shelling (and flag) vectors of $(n-1)$-graphs on an unlimited number
of vertices have a span with unlimited dimension.  Thus, $\ftilde G$ as a
sum of words must be constructed out of an infinite alphabet of letters. 
This is most likely extremely inconvenient.  The passage from $\ftilde
L_l$ (or $fL_l$) to the vertex contribution $f'L_l$ reduces this
alphabet.  Loosely speaking, it corresponds to taking the tail part of the
shelling (or flag) vector of the link.  It is a truncation.

Clearly, there are many possible ways of truncating $\ftilde L_l$, so that
the result is finite dimensional.  Each such rule gives rise to a possibly
different nullspace or abstract flag vector.  Finding the correct rule is
an important problem.  The \emph{disjoint pair of options} rule presented
in \S2 ensures that a property useful in the study of $2$-graphs continues
to hold.  (It is something like taking a linear approximation.) For each
$i$, the properties of $(i-1)$-graphs must be well understood, before the
flag vector of an $i$-graph can be given an explicit form.

For $G$ a $3$-graph on $N$ vertices, its flag vector $fG$ can be written
as a sum of length $N$ words in $a$, $b$ and $c$.  For each shelling, as
each vertex is removed, record the following.  First, an `$a$' for the
vertex itself.  Second, a `$b$' for each cell that is removed.  Third, for
each `pair of cells meeting along an edge' that is removed, record a
`$c$'.  This follows from properties of $1$-graphs. As noted, until
$3$-graphs are investigated, the flag vector of $4$-graphs cannot be given
an explicit form.

Topology also studies manifolds that possess additional structure. 
Suppose, for example, that $M$ is an \emph{oriented} topological manifold.
In that case, a triangulation of $M$ will produce an \emph{oriented
hypergraph} $G$.  This means that each cell is given an orientation (a
means of attaching a sign to each ordering of its vertices, that is
invariant under even permutations).  As before, one can define a shelling
vector $\ftilde G$ for oriented $i$-graphs.  This is because at each stage
in the shelling, the change can be expressed as the cone on something
simpler, namely an oriented $(i-1)$-graph.  In the oriented case, however,
the induction will be based on the oriented $0$-graphs.  These can be
denoted by $[a]$, $[b_+]$ and $[b_-]$.  (Here $[b_+]$ is $[b]$, which has
$\emptyset$ as its only cell, with `$+$' attached as sign to the only
ordering of its vertices.)

One can now define the \emph{flag vector of an oriented hypergraph}.  As
before, use the recursive formula (\ref{eqn:tensor-f}).  As before, the
link contribution $f'L_l$ is $fL_l$ modulo disjoint pairs of options. 
(When, as here, a set of vertices is the support for several different
cells, a basic optional cell is where one has either some cell ($+$) or no
cell ($-$) supported on these vertices.  Other options can be built up
using addition.)

It should now be clear that shelling and flag vectors exist more generally.
Let $\calX$ denote any class of objects that can be described as a union
of cells (possibly with additional structure) on a finite set of vertices.
Suppose also that the change that occurs at each stage of a shelling can
be described as the cone on something.  That something should be of some
type $\calX'$, for which a shelling vector has already been defined.  In
this situation, $\calX$ will also have a shelling vector.  If, in
addition, $\calX'$ has a flag vector, then the disjoint pair of options
rule will define a flag vector for $\calX$ also.

Thus, for example, each triangulated manifold with boundary will have both
a shelling and flag vector.  The $0$-objects that are the basis for the
induction can be represented as $[a]$, $[b]$ and $[\bdot]$.  Both $[a]$
and $[b]$ are as before.  However, $[\bdot]$ when coned twice will produce
first a special sort of $1$-object, and then a $2$-object (representing a
$1$-cell) which has the second apex as its `boundary'.  (Note that the
concept of a $0$-manifold with boundary is somewhat bizarre.)  As before,
impose `$+$' and `$-$' subscripts on the `$b$' terms to obtain the
oriented form.

Suppose $M$ is a triangulated \emph{differential} manifold.  In that case,
an oriented matroid can be used to provide a combinatorial record of the
differential structure.  In this case also, surely it is true that
shelling and flag vectors can be defined.  This the author has not
investigated.  The Pontrjagin numbers, one would wish to be linear
functions of the flag vector.  By virtue of \cite{bib:IMG+RM}, this
problem may be accessible.

It is not necessary that the classes $\calX$ and $\calX'$ have the same
general type.  Manifolds with boundary provide a mild example of this. 
Here is another.  Let $G$ be a finite group, considered as the ternary
relation $xy=z$ (or, if one prefers, $xyz=1$).  As before, a ternary
relation can be shelled.  The change at each stage of the shelling is once
again the cone on something, but what that something is requires some
thought.  Indeed, for a ternary relation the very concept of a cell
requires some thought.  The following appears to be correct.  Each
solution $S$ to the equation $xy=z$ (or triple in the ternary relation)
has a \emph{support}, which are the vertices that appear in $S$.  A
\emph{cell} will consist of a (non-empty) set of solutions $S$, each of
which has the same support.  A cell is \emph{prime} if it contains a
single solution.  Thus, for example, if $a$ and $b$ commute and $ab=c$,
then $\{\,ab=c, ba=c\,\}$ is a cell, supported on $\{a,b,c\}$.  For an
arbitrary ternary relation, there are $3!=6$ prime cells supported on a
triple $\{a,b,c\}$, and so $2^6-1=63$ possible cells altogether. 
Supported on a pair $\{a,b\}$ there are again $6$ prime cells, and so $63$
possible cells.  Supported on a singleton $\{a\}$ there is of course only
one cell, namely $\{\,aa=a\,\}$.

Consider now the shelling of a ternary relation $R$.  A prime cell such as
$\{(a,b,c)\}$ can be thought of as the cone on the prime cell
$\{(?,b,c)\}$, where $a$ is the apex of the cone.  This then is an example
of a prime cell, supported on the pair $\{b,c\}$ of vertices, for the link
(at $a$) of a ternary relation.  Once again, there will be six possible
prime cells for the link, supported on $\{b,c\}$.  Similarly, there will
be six prime cells for the link, supported on $\{b\}$.  To deal with
$\{(a,a,a)\}$ one will, for the link, need the concept of a cell supported
on an empty set of vertices.  Such will, in  a shelling, be removed at
step zero, before any other cell.  Proceeding in this way one will obtain
shelling and flag vectors for ternary relations.

(The prime $0$-objects, that are the basis for the induction, are symbols
such as $[123]$ or $[223]$, with the property that
\[
\begin{array}{ll}
        C_a\,C_b\,C_c\, [123]\,=\, \{(c,b,a)\} \> ;
    &   C_a\,C_b\,C_c\, [312]\,=\, \{(a,c,b)\} \> ;  \\
        C_a\,C_b\,C_c\, [232]\,=\, \{(b,a,b)\} \> ;
    &   C_a\,C_b\,C_c\, [333]\,=\, \{(a,a,a)\} \> ;
\end{array}
\]
where $C_x$ is the operator that forms a cone, and labels its apex `$x$'. 
Applied to a sequence $[\ldots]$ of numbers and vertices, it replaces
`$1$' by `$x$', and lowers the remaining numbers by one.  Certain symbols,
such as $[121]$, are not needed, and will not be used.  Thus,
$C_c[232]=[121]$; $C_b[121]=[b1b]$; and $C_a[b1b]=[bab]$, which is the
prime cell as above. As usual, use the disjoint pair of optional (and not
necessarily prime) cells rule to define the link contribution $f'L_l$.)

The resulting flag vector will clearly be quite complex.  This is probably
in its favour, for a group is also quite complex, and so a simple flag
vector would not be able to take a good grasp of its structure.  The
author has not investigated this matter.

It should now be clear that one would expect more or less anything that
can be shelled to have a flag vector.  Suppose the realisation $|G|$ of
$G$ is a (possibly oriented and/or differential) manifold.  Now consider
the question: \emph{what part of the flag vector $fG$ is a topological
invariant of $M=|G|$?}  To understand how barycentric subdivision will
change $fG$ is an obvious starting point, for topological (differential)
invariants will vanish on such changes to $fG$.  Morally, although perhaps
not in fact, this is all that is required, to produce the class of
invariants that corresponds to the given $f$-vector.

Now note that a `$G$  with a barycentric subdivision' is an object $\Ghat$
that can be shelled.  Presumably, there is a corresponding flag vector
$f\Ghat$, from which both $fG$ and $fG'$ (the result of the subdivision)
can be computed, and thus the change $f(\partial\Ghat) = fG'-fG$ in $fG$
also.

Objects $G$ whose realisation $|G|$ is a manifold are rather special, and
in general their flag vectors $fG$ will span a subspace of all flag
vectors.  (The precise determination of the subspace is an important
problem, which may lead to crossing-change rules.  The first step is to
produce a spanning set for the manifold nullspace $\calM_0$, defined in
\S1.)  The same applies to $\Ghat$, namely `$G$ with a barycentric
subdivision'.  If the properties of the flag vector allow these subspaces
to be determined then, at least morally, $fG$ modulo all possible
$f(\partial\Ghat)$ is an invariant of $M=|G|$.

Before concluding, here are some comments that do not conveniently belong
elsewhere.  First, the domain of the theory.  By using a suitable notion
of a cell, one can triangulate and shell any of: a submanifold $N\subset
M$; a singular stratified space $X$; the realisation of a simplicial
complex.  Thus, for example, one can obtain a class of invariants for
knots in $\bfR^3$.  For knots in a $3$-manifold $M$ more work is required,
for here homeomorphism for $N\subset M$ will in general be a weaker
concept than isotopy.  To allow a $4$-manifold to have quantum topology,
but not singular algebraic surfaces is, in the author's view, not
reasonable.  Infinite objects can be accomodated, provided some factor is
inserted into (\ref{eqn:tensor-f}), so as to ensure convergence.  In
addition, only `finite partial shellings' should be used.

Flag vectors are very important in the study of convex polytopes (and the
associated algebraic varieties).  Here, there are four steps.  First, the
definition of the flag vector.  Second, the linear equations on the flag
vector.  Third, the linear inequalities and the associated homology
theory.  Fourth, the pseudopower (usually non-linear) inequalities and the
`ring-like' structure.

(For polytopes, the first step is all but trivial.  The flag vectors of
dimension $n$ polytopes span a space, whose dimension is
\cite{bib:BB,bib:MVIC} the $(n+1)$st Fibonacci number.  Middle perversity
intersection homology and, the author believes, the local-global extension
\cite{bib:LGIH} is the associated homology theory, while \cite{bib:IHRS}
may provide the `ring-like' structure.  For simple polytopes, the entire
theory has a satisfactory form.  This possible analogy between convex
polytopes and $i$-graphs has been, for the author, an important motive.)

For $i$-graphs, this paper has considered only the first step.  For
topology, there is a fifth step, which is to determine which linear
functions $v$ of the $f$-vector are topological (differential)
invariants.  As with the Kontsevich knot integral, part of $f$ will vary
with the representation of $M=|G|$.  As noted for $2$-graphs, if there is
a suitable homology theory associated to the flag vector, it will produce
subtle geometric inequalities.

The quantum $3$-manifold invariants \cite{bib:LMO,bib:EW-QFTJP} can be
expressed as a formal sum of trigraphs (with vertex orientations) modulo
certain relations.  One would like to be able to express this theory in
terms of flag vectors.  A first step is this.  If $G$ is an oriented
$4$-graph (perhaps a triangulation of a $3$-manifold), one can produce
trigraphs from $G$ in the following way. The trigraph will be realised in
$|G|$. Each edge will pass through the centroid of a facet of a cell. 
Within each cell (or rather its realisation) one might have nothing, a
segment linking two facet centroids, a pair of such segments, or a
trivalent vertex (oriented by the cell) connecting three facet centroids.
Determining the coefficients is of course another matter.  The relations
should correspond, of course, to barycentric subdivision.

For an oriented compact and connected $2$-manifold $M$, the Euler
characteristic $\chi(M)$ is the only topological invariant, and it can be
computed by integrating the curvature.  If however $M$ is not connected,
and such cases should be considered, it has a $\chi$ for each connected
component. Now triangulate $M$, to solve $|G|=M$.  The Euler
characteristic is on $G$ a sum of local contributions.  It may be that
`facet centroid to adjacent facet centroid' loops will allow these local
contributions to be glued together in a quantum (or flag) way, so that
they contribute `only to the same component'.  Looking the other way,
inside a triangulated $4$-manifold one can see embedded $3$-graphs
(singular surfaces).

There are some points of contact between the flag vector and the concept
of a topological quantum field theory (TQFT) \cite[Chap.~2]{bib:MFA-GPK}. 
Briefly, suppose a rule $v=v(M)$ that assigns a vector to every (oriented)
$n$-manifold has been given.  Now, for every (oriented) $(n-1)$-manifold
$B$ define a space $Z(B)$ as follows.  Take all formal sums of manifolds
with boundary $B$, modulo the following relation.  A formal sum
\[
    M_1 = \lambda_1 M_{1,1} + \ldots + \lambda_r M_{1,r}
\]
is treated as zero if
\[
    v(M_1 \cup_BM_2) = \lambda_1 v(M_{1,1}\cup_BM_2) + \ldots + 
                        \lambda_r v(M_{1,r}\cup_BM_2)
\]
is zero for every manifold $M_2$ (or formal sum of such) with boundary
$B^*$ (the boundary $B$ of $M_1$, but with the reverse orientation).  One
can then ask of $v=v(M)$, whether its associated $Z(B)$ spaces satisfy
properties such as the multiplicative and associative axioms of TQFT.  A
similar construction can be defined, and thus similar questions asked, for
the flag vector of $i$-graphs, and suitable functions thereof.

(The definition of $G=G_1\cup_BG_2$ requires some thought. The following
may work the best.  Suppose the vertex set of $G$ is partitioned into
$V_1$ and $V_2$.  Now define $G_i$ to be the subgraph of $G$, whose cells
have support a subset of $V_i$.  Define $B$ to be the $i$-graph consisting
of the remaining cells. This corresponds to a `collared' form of the
manifold construction.  Note that the cells of $G_1$ are disjoint from
those of $G_2$.)

Finally, to conclude here is a concise summary of the whole paper.  Any
finite object $G$ built out of `cells' can be shelled.  By induction, this
gives rise to a shelling vector $\ftilde G$.  The same can be done for
some infinite $G$ also.  The $0$-objects, that are the basis of the
induction, depend on the type of object that $G$ is.  For each rule
$L_l\mapsto f'L_l$ defining a link contribution $f'L_l$, the recursive
formula
\[
    fG = \sum \nolimits _{\rm shellings} \>\>
            f'L_1 \otimes \ldots \otimes f'L_N
\]
defines a flag vector $fG$.  This paper, without investigating the matter
closely, proposes that $f'L_l$ should be $fL_l$ modulo disjoint pairs of
optional cells.  In at least some cases the resulting spaces (of flag
vectors and link vectors) can be given a pleasant and explicit
description.  Indeed, to do this, without discarding vital information, is
perhaps the main purpose of the rule for the link contribution $f'L_l$. 
When $G$ has a topological (or differential) realisation $|G|$, the
topological (differential) invariants that are a linear function of the
flag vector $fG$ can be studied, particularly via barycentric subdivision. 
Whether or not these invariants are related to those of quantum topology
requires further investigation.

\end{document}